\documentclass[12pt]{article}
\usepackage{amsmath}
\usepackage{amssymb}
\usepackage{theorem}

\newcommand{\eop}{\bigstar}

\newcommand{\card}[1]{{\vert #1 \vert} }

\newcommand{\otp}[1]{\hbox{otp($#1$)}}

\newcommand{\Dom}{{\rm Dom}}

\newcommand{\llg}{{l\rm g}}

\newcommand{\Rang}{{\rm Rang}}

\newcommand{\cf}{{\rm cf}}

\newcommand{\initial}{\vartriangleleft}

\newenvironment{Proof}{\noindent{\bf Proof.}}{\par\bigskip} 

{\theorembodyfont{\rmfamily}\newtheorem{THEOREM}{Theorem}[section]}

{\theorembodyfont{\rmfamily}\newtheorem{Conclusion}[THEOREM]{Conclusion}}

{\theorembodyfont{\rmfamily}\newtheorem{LEMMA}[THEOREM]{Lemma}}

{\theorembodyfont{\rmfamily}\newtheorem{Main Theorem}[THEOREM]{Main Theorem}}
\newenvironment{main Theorem}{\begin{Main Theorem}} 
{\end{Main Theorem}}

{\theorembodyfont{\rmfamily}\newtheorem{Preservation Lemma}[THEOREM]{Preservation Lemma}}
\newenvironment{preservation lemma}{\begin{Preservation Lemma}}
{\end{Preservation Lemma}}

{\theorembodyfont{\rmfamily}\newtheorem{Construction Lemma}[THEOREM]{Construction Lemma}
\newenvironment{construction lemma}{\begin{Construction Lemma}}
{\end{Construction Lemma}}

{\theorembodyfont{\rmfamily}\newtheorem{Theorem}[THEOREM]{Theorem}}

{\theorembodyfont{\rmfamily}\newtheorem{Definition}[THEOREM]{Definition}}

{\theorembodyfont{\rmfamily}\newtheorem{Conventions}[THEOREM]{Conventions}}

{\theorembodyfont{\rmfamily}\newtheorem{Main Definition}[THEOREM]{Main Definition}}
\newenvironment{main definition}{\begin{Main Definition}}
{\end{Main Definition}}

{\theorembodyfont{\rmfamily}\newtheorem{Lemma}[THEOREM]{Lemma}}

{\theorembodyfont{\rmfamily}\newtheorem{Notation}[THEOREM]{Notation}}

{\theorembodyfont{\rmfamily}\newtheorem{Convention}[THEOREM]{Convention}}

{\theorembodyfont{\rmfamily}\newtheorem{Note}[THEOREM]{Note}}

{\theorembodyfont{\rmfamily}\newtheorem{Observation}[THEOREM]{Observation}}

{\theorembodyfont{\rmfamily}\newtheorem{Remark}[THEOREM]{Remark}}

{\theorembodyfont{\rmfamily}\newtheorem{Main Fact}[THEOREM]{Main Fact}}
\newenvironment{main Fact}{\begin{Main Fact}}{\end{Main Fact}}

{\theorembodyfont{\rmfamily}\newtheorem{Fact}[THEOREM]{Fact}}

{\theorembodyfont{\rmfamily}\newtheorem{Subfact}[THEOREM]{Subfact}}

{\theorembodyfont{\rmfamily}\newtheorem{Claim}[THEOREM]{Claim}}

{\theorembodyfont{\rmfamily}\newtheorem{Main Claim}[THEOREM]{Main Claim}}
\newenvironment{main claim}{\begin{Main Claim}}{\end{Main Claim}}

{\theorembodyfont{\rmfamily}\newtheorem{Corrolary}[THEOREM]{Corrolary}}

{\theorembodyfont{\rmfamily}\newtheorem{Subclaim}[THEOREM]{Subclaim}}

{\theorembodyfont{\rmfamily}\newtheorem{Corollary}[THEOREM]{Corollary}}

{\theorembodyfont{\rmfamily}\newtheorem{Example}[THEOREM]{Example}}

{\theorembodyfont{\rmfamily}\newtheorem{Proposition}[THEOREM]{Proposition}}

{\theorembodyfont{\rmfamily}\newtheorem{Discussion}[THEOREM]{Discussion}}

\newenvironment{Proof of the Subfact}
{\noindent{\bf Proof of the Subfact.}}{\par\bigskip}

\newenvironment{Proof of the Theorem}
{\noindent{\bf Proof of the Theorem.}}{\par\bigskip}

\newenvironment{Proof of the Conclusion}
{\noindent{\bf Proof of the Conclusion.}}{\par\bigskip}

\newenvironment{Proof of the Observation}
{\noindent{\bf Proof of the Observation.}}{\par\bigskip}

\newenvironment{Proof of the Fact}
{\noindent{\bf Proof of the Fact.}}{\par\bigskip}

\newenvironment{Proof of the Lemma}
{\noindent{\bf Proof of the Lemma.}}{\par\bigskip}

\newenvironment{Proof of the Claim}
{\noindent{\bf Proof of the Claim.}}{\par\bigskip}

\newenvironment{Proof of the Subclaim}
{\noindent{\bf Proof of the Subclaim.}}{\par\medskip}

\newenvironment{Proof of the Main Claim}
{\noindent{\bf Proof of the Main Claim.}}{\par\bigskip}

\newcommand{\elementary}{\prec}

\newcommand{\Bbf}{\Bbb}

\newcommand{\into}{\rightarrow}

\newcommand{\rest}{\upharpoonright}  

\newcommand{\nacc}{\mathop{\rm nacc}}

\newcommand{\cl}{\mathop{\rm cl}}
\newcommand{\acl}{\mathop{\rm acl}}
\newcommand{\dcl}{\mathop{\rm dcl}}

\newcommand{\tp}{\mathop{\rm tp}}

\newcommand{\deq}{\buildrel{\rm def}\over =}

\newcommand{\FF}{{\cal F}}

\newcommand{\PP}{{\cal P}}

\title{On properties of theories which preclude the existence of universal
models}
\author{Mirna D\v zamonja\\
School of Mathematics\\
University of East Anglia\\
Norwich, NR4 7TJ, UK
\\
\scriptsize{M.Dzamonja@uea.ac.uk}\\
\scriptsize{http://www.mth.uea.ac.uk/people/md.html}\and
Saharon Shelah\\
Mathematics Department\\
Hebrew University of Jerusalem\\
91904 Givat Ram, Israel\\
and\\
Mathematics Department
Rutgers University\\
New Brunswick, New Jersey, USA\\
\scriptsize{shelah@sunset.huji.ac.il}\\
\scriptsize{http://www.math.rutgers.edu/$\thicksim$shelarch}
}

\begin{document}

\baselineskip=16pt
\binoppenalty=10000
\relpenalty=10000
\raggedbottom

\maketitle

\begin{abstract} In this paper we investigate some properties of
first order theories which prevent them from having universal
models under
certain cardinal arithmetic assumptions. 
Our results 
give a new syntactical condition, oak property, which is a sufficient condition
for a theory not to have universal
models in cardinality $\lambda$ when certain cardinal arithmetic assumptions
implying the failure of $GCH$ (and close to the failure of
$SCH$) hold.

\footnote{
This publication is numbered 710 in the list of publications
of Saharon Shelah. 
The authors thank the United-States Israel Binational Science
Foundation and NSF for their support
during the preparation of this paper, and Mirna D\v zamonja thanks
the Academic Study Group for their support during the summer of 1999
and Leverhulme Trust for their grant number number F/00204B.


}
\end{abstract} 

\section{Introduction} The existence of a universal model of
a theory has been the object of a continuous interest to
specialists in various disciplines of mathematics,
see for example \cite{AB}, \cite{Komjath}. We approach
this problem from the point of view of model theory, more
specifically, classification theory, and
we concentrate on first order theories. In a series of papers,
Kojman-Shelah \cite{KjSh 409} (see there also for earlier
references), \cite{KjSh 447}, Kojman \cite{Kj},
Shelah \cite{Sh 457},
\cite{Sh 500}, D\v zamonja-Shelah \cite{DjSh 614}, the thesis 
claiming the connection between the complexity of a theory
and its amenability to the existence of universal models,
has been pursued. As it follows from the classical results in model
theory (see \cite{ChKe}) that
if $GCH$ holds, then every countable first order
theory admits a universal model in every uncountable cardinal,
the question we need to ask is what happens when $GCH$ fails.
It is usually ``easy" to force a situation in which there are no universal
models (by adding Cohen subsets), however assuming that $GCH$ fails and
allowing ourselves a vague use of the
words ``many"
and ``often", we can distinguish between those theories which
for many cardinals do not have a universal model in that
cardinal {\em whenever} $GCH$ fails, and those for which it is
possible to
construct a model of set theory in which $GCH$ fails, yet our
theory has a universal model in the cardinality under consideration.
This division would suggest that  the theories first described,
let us call them for the sake of this introduction amenable,
are of higher complexity than the latter ones. 

In his paper \cite{Sh 500}, S. Shelah introduced a hierarchy
of complexity for first order theories, and showed that
past a certain level on that hierarchy, the inherent properties
of any theory on that level, will preclude the existence of
universal models in most cardinalities. The details of this
hierarchy are described in the following Definition \ref{nsops},
and what S. Shelah proved in \cite{Sh 500}, is that
$SOP_4$ implies high non-amenability. Here we
show that this bound is not sharp, by defining a property
of theories which is present in some $NSOP_4$ theories
(meaning, not $SOP_4$),
yet it precludes the existence of a universal model under
certain cardinal arithmetic. 
This property is called the oak property, as its
prototype is the model completion of
Th$(M_{\lambda,\kappa})$, a theory connected to that of the tree
${}^{\kappa\ge}\lambda$ (for details see Example \ref{prototype}).
The oak property cannot be made
a part of the $SOP_n$ hierarchy,
as we exhibit a theory which has oak, and is $NSOP_3$, while the model
completion of the theory of triangle free graphs is an example of a $SOP_3$
theory which does not satisfy the oak property.
Our research is a continuation of section \S1 of \cite{Sh 457}, where
the universality spectrum of the theory
$T_{\rm feq}^\ast$ of
infinitely many indexed independent equivalence relations is investigated, and
it is proved that under cardinals arithmetic assumptions like the ones in our
Theorem \ref{main}, $T_{\rm feq}^\ast$ does not have universal models.
We show that $T_{\rm feq}^\ast$ has the oak property,
and in fact exhibit a close connection between $T^\ast_{\rm feq}$
and Th$(M_{\lambda,\kappa})$.


We commence by giving some background notions which will be
used in the main sections of the paper.
First, several classical definitions of model
theory.

\begin{Convention} A theory in this paper
means a first order complete theory, unless otherwise stated. Such an object is
usually denoted by $T$.
\end{Convention}

\begin{Notation} Given
a theory $T$, we let $\mathfrak C={\mathfrak C}_T$
stand for ``the monster model", i.e. a saturated enough model of
$T$. As is usual, we assume without loss of
generality that all our discussion takes place inside some such model,
so all expressions to the extent ``there is", ``exists" and ``$\models$"
are to be relativised to this model, all models are $\elementary {\mathfrak C}$,
and all subsets of $\frak C$ we mention have size less than the saturation
number of $\frak C$.
We let
$\bar{\kappa}=\bar{\kappa}({\mathfrak C}_T)$ be the size of $\mathfrak C$, so
this cardinal is larger than any other cardinal mentioned in connection with
$T$. 
\end{Notation}

\begin{Definition} (1) The tuple $\bar{b}$ is {\em defined by}
$\varphi(\bar{x};\bar{a})$ if $\varphi(\mathfrak{C};\bar{a})=\{\bar{b}\}$,
i.e. if $\bar{b}$ is the unique $\bar{x}$ which realizes $\varphi(\bar{x};
\bar{a})$.
It is defined by the type $p$ if $\bar{b}$ is the unique tuple which
realizes $p$. It is definable over $A$ if $\tp(\bar{b},A)$ defines it.

(2) The formula $\varphi(\bar{x};\bar{a})$ is {\em algebraic} if
$\varphi(\mathfrak{C};\bar{a})$ is finite. The type $p$ is
algebraic if it is realized by finitely many tuples only. The
tuple $\bar{b}$ is algebraic over $A$ if $\tp(\bar{b},A)$ is.

(3) The {\em definable closure} of $A$ is
\[
\dcl(A)\deq\{b:\,b\mbox{ is definable over }A\}.
\]

(4) The {\em algebraic closure} of $A$ is
\[
\acl(A)\deq\{b:\,b\mbox{ is algebraic over }A\}.
\]

(5) If $A=\acl(A)$, we say that $A$ is {\em algebraically closed}.
When $\dcl(A)$ and $\acl(A)$ coincide, then $\cl(A)$ denotes their common
value. \end{Definition}

\begin{Definition} (1) For a theory $T$ and a cardinal $\lambda$,
models $\{M_i:\,i<i^\ast\}$ of $T$, each of size $\lambda$,
are {\em jointly universal} iff for every $N$ a model of $T$ of size
$\lambda$, there is an $i<i^\ast$ and an isomorphic embedding
of $N$ into $M_i$.

{\noindent (2)} For $T$ and $\lambda$ as above, 
\begin{equation*}
\begin{split}
{\rm univ}(T,\lambda)\deq\min\{\card{\FF}:\,\FF\mbox{ }
&\mbox{is a family of
jointly}\\ &\mbox{universal models of }T\mbox{ of size }\lambda\}.\\
\end{split}
\end{equation*}
{\scriptsize{(so univ$(T,\lambda)=1$ iff there us a universal model of $T$ 
of size $\lambda$.)}}
\end{Definition}
The following is the main definition of S. Shelah's \cite{Sh 500}.

\begin{Definition}\label{nsops} (Shelah, \cite{Sh 500}) Let $n\ge 3$.
\begin{description}
\item{(1)} A formula $\varphi(\bar{x},\bar{y})$
is said to exemplify the $n$-{\em strong order
property},
$SOP_n$ if $\llg(\bar{x})=\llg(\bar{y})$, and there are
$\bar{a}_k$  for $k<\omega$, each of length ${\llg(\bar{x})}$
such that

\begin{description}
\item{(a)} $\models\varphi[\bar{a}_k,\bar{a}_m]$ for $k<m<\omega$,

\item{(b)} $\models \neg(\exists\bar{x}_{0},\ldots,\bar{x}_{n-1})[\bigwedge\{
\varphi(\bar{x}_{l},\bar{x}_{k}):\,l,k<n\mbox{ and }k=l+1\mod n\}].$

\end{description}

$T$ has $SOP_n$
if there is a formula $\varphi(\bar{x},\bar{y})$ exemplifying this.

\item{(2)} $SOP_{\le n}$ is defined similarly, except that in (b) we replace
``$n$" by each ``$m\le n$".

\item{(3)} $NSOP_n$ stands for the negation of $SOP_n$.

\end{description}
\end{Definition} 

\begin{Note} Using a compactness argument and Ramsey theorem, one can prove
that
if $T$ is a theory with $SOP_n$ and $\varphi(\bar{x},\bar{y})$,
and $\langle \bar{a}_n:\,n<\omega\rangle$ exemplify it, without loss
of generality $\langle \bar{a}_n:\,n<\omega\rangle$ is an indiscernible
sequence. See \cite{Sh -c}, or \cite{Rami} for examples of such arguments.
\end{Note}

\begin{Example} The model
completion of the theory of triangle-free graphs is a
prototypical example of a $SOP_3$ theory,
with the formula $\varphi(x,y)$ just stating that $x$ and $y$ are connected.
It can be shown that this theory is $NSOP_4$, see \cite{Sh 500}.
\end{Example}

The following fact indicates that $SOP_n (3\le n<\omega)$ form a
hierarchy, and the thesis is that this hierarchy is reflected in
the complexity of the behavior of the relevant theories under natural
constructions in model theory.

\begin{Fact}\label{decreasing}(Shelah, \cite{Sh 500}, \S2) $SOP_{n+1}\implies SOP_n$.
\end{Fact}

\section{An $NSOP_3$ theory without universals}\label{NSOP3}

\begin{Definition}\label{T} (1)
Let $T_0$ be the following theory in the language 
\[
\{Q_0,Q_1,Q_2, F_0, F_1, F_2, F_3\}:
\]

\begin{description}
\item{(i)}
$Q_0,Q_1,Q_2$ are unary predicates which form a partition of the universe,
\item{(ii)}
$F_0$ is a partial function from $Q_1$ to $Q_0$,
\item{(iii)}
$F_1$ is a partial two-place function from $Q_2\times Q_0$ to
$Q_1$.
\item{(iv)} $F_2$ is a partial function from $Q_0$ to $Q_2$,
\item{(v)} $F_3$ is a partial function from $Q_2$ to $Q_0$,
\item{(vi)}
$F_0(F_1(z,x))=x$ for all $(z,x)\in \Dom(F_1)$, and
\item{(vii)} $F_3(F_2(x))=x$ for all $x$.
\end{description}

{\noindent (2)} Let $T_0^+$ be like $T_0$, but with the requirement that
$F_0, F_1, F_2$ and $F_3$ are total functions.
\end{Definition}

\begin{Remark} It is to be noted that the above definition
of $T_0$ uses partial, 
rather than the more usual, full function symbols.
Using partial functions, we have to be careful when we
speak about submodels, 
where we have a choice of deciding whether statements of the form ``$F_l(x)$ is
undefined" are preserved in the larger model. We choose to
request that the fact that $F_l$ is undefined at a certain entry, is
not necessarily preserved in the larger
model. Functions $F_2$ and $F_3$ are ``dummies" whose sole purpose is to assure
that models of $T^+_0$ are non-trivial, while keeping $T_0^+$ a universal
theory (which is useful when discussing the model completion).
 Also note that neither $T_0$ nor $T_0^+$ is complete, but
every model $M$ of $T_0$ in which $Q_0^M,Q_2^M\neq\emptyset$ and $F_0$
and $F_3$ are onto,
can be extended to a model of $T_0^+$ with the same universe (Claim \ref{amalg}
(2)), and every model of $T_0$ is a submodel of a model of $T^+_0$ (Claim
\ref{amalg}(4)).
$T_0^+$ has a complete model completion (Claim \ref{mcompletion}). This model completion is
the main theory we shall work with and,
as we shall show, it has the oak property  (Claim \ref{nasa}) and is $NSOP{}_4$
(Claim \ref{nsop}).
\end{Remark}

\begin{Example}\label{prototype} An example which we take as the prototype 
of a model of $T_0^+$,
is a model $M=M_{\lambda,\kappa}$ obtained when for given
infinite cardinals $\kappa,\lambda$, we take
$Q_0^M$ to be $\kappa$, $Q_1^M$ to be ${}^{\kappa>}\lambda$,
and $Q_2^M={}^\kappa\lambda$. We let $F_0(\eta)$ be the length of $\eta$ for
$\eta\in Q_1$, and let $F_1(\nu,\alpha)=\nu\rest\alpha$. Let 
$F_3$ be any surjective function from $Q_2^M$ onto $Q_0^M$, and
for $\alpha<\kappa$ let $F_2(\alpha)=\nu_\alpha$ for any $\nu_\alpha$
such that $F_3(\nu_\alpha)=\alpha$.
\end{Example}

\begin{Claim}\label{amalg}
\begin{description}
\item{(1)} If $M$ is a model of $T_0^+$, then $Q_0^M$, $Q_1^M$ and
$Q_2^M$ are all non-empty, and $F_0^M$ and $F_3^M$ are onto.

\item{(2)} Every model $M$ of $T_0$
in which $Q_0^M\neq\emptyset$ and $Q_2^M\neq\emptyset$, while $F_0$
and $F_3$ are onto, can be
extended to a model of $T_0^+$ with the same universe (and every model of
$T_0^+$ is a model of $T_0$).

\item{(3)} There are models $M$ of $T_0$ with $Q_0^M\neq\emptyset$ and $Q_2^M\neq\emptyset$
and $F_3^M$ onto,
which cannot be extended to a model of $T_0^+$ with the same universe.

\item{(4)} Every model of $T_0$ is a submodel of a model of $T_0^+$.

\item{(5)} $T_0^+$ has the amalgamation
property and the joint amalgamation property
$JEP$. 

\item{(6)} If $M\models T_0$ and $A\subseteq M$ is finite, \underline{then}
the closure of $B$ of $A$ under $F_0, F_1, F_2$ and $F_3$ is finite
(in fact $\card{B}\le 12\card{A}^2+ 8\card{A}$), moreover:
\begin{description}
\item{(a)} $B\cap Q_2^M=(A\cap Q_2^M)\cup\{F_2(a):\,a\in A\cap Q_0^M\}$,
\item{(b)} $B\cap Q_0^M=(A\cap Q_0^M)\cup\{F_0(b):\,b\in A\cap Q_1^M\}
\cup\{F_3(c):\,c\in A\cap Q_2^M\}$ and
\item{(c)} $B\cap Q_1^M=
(A\cap Q_1^M)\cup\{F_1(c,a):\,c\in B\cap Q_2^M\,\,\&\,\,
a\in B\cap Q_0^M \}.$
\end{description}
In this case, $B\models T_0$ and if $M\models T^+_0$, then $B\models T_0^+$.
\end{description}
\end{Claim}

\begin{Proof of the Claim}\begin{description}
\item{(1)} As $M$ is a model we have that $M\neq\emptyset$, so at least
one among $Q_0^M, Q_1^M, Q_2^M$ is not empty.

If $Q_0^M\neq\emptyset$, then $F_2$ guarantees that $Q_2^M\neq\emptyset$, so
$Q_1^M\neq\emptyset$ because of $F_1$.
If $Q_1^M\neq\emptyset$, then
$Q_0^M\neq\emptyset$ because of $F_1$. Finally, if $Q_2^M\neq\emptyset$, then
$Q_0^M\neq\emptyset$ because of $F_3$, and we can again argue as above.

If $a\in Q_0^M$, let $b\in Q_2^M$ be arbitrary. Then $F_1(a,b)\in Q_1^M$
and $F_0(F_1(a,b))=a$. Hence, $F_0$ is onto. Also, $F_3(F_2(a))=a$, so $F_3^M$
is onto. 

\item{(2)} Let $M\models T_0$ and $Q_0^M, Q_2^M\neq\emptyset$. For $x\in Q_0^M$
and $z\in Q_2^M$ let $F_1(z,x)=y$ for any $y\in Q_1^M$ such that $F_0(y)=x$,
which exists as $F_0^M$ is already onto. For $x\in Q_0^M$ for which
$F_2(x)$ is not already defined, let $F_2(x)=z$ for any $z$ such that
$F_3(z)=x$, which exists as $F_3^M$ is onto. Finally,
extend $F_0$ and $F_3$ to be total. The described model is a model of $T_0^+$.

\item{(3)} Let $\kappa_1<\kappa_2<\lambda$ and let $Q_0^M=\kappa_2$,
$Q_1^M={}^{\kappa_1>}\lambda$, while $Q_2^M={}^{\kappa_1}\lambda$. For $\alpha
<\kappa_2$ let $F_2(\alpha)$ be the function
in ${}^{\kappa_1}\lambda$ which is constantly $\alpha$, and for $\nu\in 
{}^{\kappa_1}\lambda$ let $F_3(\nu)=\min(\Rang(\nu_2))$ if this value is
$<\kappa_2$, and $0$ otherwise. Also, let $F_0(\eta)=\llg(\eta)$ and $F_1
(\nu,\alpha)=\nu\rest\alpha$ be defined for $\nu\in{}^{\kappa_1}\lambda$
and $\alpha<\kappa_1$.

This is a model of $T_0$, but not of $T_0^+$ because $F_1$ is not total.
If this model were to be extended to a model of $T_0^+$ with the same
universe, we would have that for every $\nu\in {}^{\kappa_1}\lambda$
\[
F_0(F_1(\nu,\kappa_1))=\kappa_1\,\,\&\,\,F_1(\nu,\kappa_1)=\eta
\]
for some $\eta\in {}^{\kappa_1>}\lambda$. As $F_0(\eta)$ is already defined,
$F_0(\eta)=\llg(\eta)<\kappa_1$, which is a contradiction.

\item{(4)} Given a model $M$ of $T_0$. First assure that $Q_0^M, Q_1^M,
Q_2^M\neq\emptyset$ by adding new elements if necessary. Then make sure that
$F_0$ and $F_3$ are total and onto, again by adding new elements if
needed. Now define $F_1(z,x)=y$ if $F_0(y)=x$, which is possible. Finally,
declare $F_2(x)=z$ for any $z$ such that $F_3(z)=x$.

\item{(5)} Suppose that $M_0, M_1$ and $M_2$ are models of $T_0^+$ with
$\card{M_1}\cap \card{M_2}=\card{M_0}$,
and $M_0\subseteq M_1, M_2$. We define $M_3$ 
as follows. Let $\card{M_3}=\card{M_1}
\bigcup \card{M_2}$, and for $m\in \{0,2,3\}$ let
$F_m^{M_3}(x)= F_m^{M_l}(x)$ if $x\in M_l$ for some $l$. This is well defined,
because $M_1$ and $M_2$ agree on $M_0$. Also, the identity $F_3(F_2(x))=x$
is satisfied in $M_3$.

For $(z,x)\in Q_2\times Q_0^M$ such that for some $l$ we have $x
\in M_l$ and $z\notin M_l$ choose $y_{z,x}\in M_l$ such that $F_0^{M_l}(y)=x$,
which is possible by part (1) of this Claim. Now we define $F_1$ by letting
or $(z,x)\in Q_2\times Q_0^M$
\begin{equation*}
F_1(z,x)=
\begin{cases}
F_1^{M_l}(z,x)& \text{if $z,x\in M_l$},\\
y_{(z,x)}& \text{otherwise}.
\end{cases}
\end{equation*}
Now it can be easily seen that $M_3$ is a model of $T^+_0$ and that both $M_1$
and $M_2$ are it submodels. This proves the amalgamation property for $T^+_0$.

To see that JEP holds, suppose that we are given two models $M_1$, $M_2$ of
$T^+_0$. We let $M$ be their disjoint union and define the functions $F_m$
for $m\in \{0,1,2,3\}$.

\item{(6)} Clearly $B$ is contained in the closure of $A$ and the size of $B$
is as claimed. It can be checked directly that $B$ is closed, using the
equations of $T_0$, and it also easily follows that $B$ is a model of $T_0$,
or of $T_0^+$ if $M$ is.

\end{description}

$\bigstar_{\ref{amalg}}$
\end{Proof of the Claim}

\begin{Claim}\label{mcompletion} $T_0^+$ has a complete model completion
$T^\ast$
which admits elimination of quantifiers, and is $\aleph_0$-categorical. In this
theory, the closure and the algebraic closure coincide.
\end{Claim}

\begin{Proof of the Claim} We can construct $T^\ast$ directly. $T^\ast$
admits elimination of quantifiers because $T_0^+$ has the amalgamation
property (\cite{ChKe} 3.5.19). It can be seen from the construction of $T^\ast$
that it is complete, or alternatively, it can be seen that $T^\ast$
has
JEP and so by \cite{ChKe} 3.5.11, it is complete.
To see that the theory is $\aleph_0$-categorical, observe that
Claim \ref{amalg}(6) implies that for every $n$
there are only finitely many $T_0$-types in $n$-variables. Then by the
Characterisation of complete $\aleph_0$-categorical theories
(\cite{ChKe} 2.3.13), $T^\ast$ is $\aleph_0$-categorical. Using the
elimination of quantifiers and the fact that all relational symbols of
the language of $T^\ast$ have infinite domains in every model of $T^\ast$,
we can see that the algebraic closure and the definable closure coincide in
$T^\ast$.
$\eop_{\ref{mcompletion}}$
\end{Proof of the Claim}

\begin{Observation}\label{insert} If $A,B\subseteq {\mathfrak C}_{T^\ast}$ are
closed
and $c\in {\rm cl}(A\cup B)\setminus A \setminus B$, \underline{then}
$c\in Q_1^{{\mathfrak C}_{T^\ast}}$.
\end{Observation}

\begin{Proof} Notice that 
\begin{equation*}
\begin{split}
{\rm cl}(A\cup B)=& A\cup B\cup
\{F_1(c,a):\,c\in (A\cup B)\cap Q_2\,\,\&\,\,
a\in (A\cup B)\cap Q_0\\
&\&\,\,\{c,a\}\nsubseteq A\,\,\&\,\,
\{c,a\}\nsubseteq B\}\\
\end{split}
\end{equation*}
by Claim \ref{amalg}(6).
\end{Proof}

\begin{Claim}\label{nsop} $T^\ast$ is $NSOP_3$, consequently $NSOP_4$.
\end{Claim} 

\begin{Proof of the Claim}
Suppose that $T^\ast$ is $SOP_3$ and let $\varphi(\bar{x},\bar{y})$,
and $\langle \bar{a}_n:\,n<\omega\rangle$ exemplify this in a model
$M$ (see Definition \ref{nsops}(1)). Without
loss of generality, each $\bar{a}_n$ is without repetition and is closed
(recall Claim \ref{amalg}(6)).
By the Ramsey theorem and compactness, we can assume
that the given sequence is a part of an indiscernible
sequence $\langle \bar{a}_k:\,k
\in {\Bbf Z}\rangle$, hence $\bar{a}_k$s form a $\Delta$-system.
Without loss of generality,
each $\bar{a}_k$ is closed under $F_0$.
Let for $k\in {\Bbf Z}$
\begin{equation*}
X_k^{<}\deq \bigcap_{m<k}{\cl}
(\bar{a}_{m}\,\hat{}\,\bar{a}_{k}),\quad
X_k^{>}\deq \bigcap_{m>k}{\cl}
(\bar{a}_{m}\,\hat{}\,\bar{a}_{k}),\quad X_k=\cl(X_k^{<}\cup X_k^{>}).
\end{equation*}
Hence ${\rm Rang}(\bar{a}_k)\subseteq X_k$, and $X_k$ is closed.
By Claim \ref{amalg}(6), there is an a priori finite bound on the size of
$X_k$,
hence by indiscernibility, we have that $\card{X_k}=n^\ast$ for some fixed
$n^\ast$ not depending on $k$.
Let $\bar{a}_k^+$ list $X_k$ with no repetition. By 
Observation \ref{insert}, Claim \ref{amalg}(6), indiscernibility and the fact
that each $\bar{a}_k$ is closed under $F_0$, we have that
\[
X_k\cap Q_0^{\mathfrak C}\subseteq \Rang(\bar{a}_k)\mbox{ and }
X_k\cap Q_2^{\mathfrak C}\subseteq \Rang(\bar{a}_k).
\]
Applying Ramsey theorem again, without loss of generality we have that
$\langle \bar{a}^+_k:\,k\in {\Bbf Z}\rangle$ are indiscernible.
Let 
\[
w^\ast_0\deq\{l:\,\bar{a}^+_{k_1}(l)=\bar{a}^+_{k_2}(l)
\mbox{ for some (equivalently all) }k_1\neq k_2\}.
\]
If $\bar{a}^+_{k_1}(l_1)=\bar{a}^+_{k_2}(l_2)$ for some $k_1\neq k_2$,
without loss of generality $k_1<k_2$, by indiscernibility. By transitivity,
using $k_1<k_2<k_3$, we get $l_1=l_2\in w^\ast_0$. Let $w_1^\ast\deq
n^\ast\setminus w^\ast_0$, and let $\bar{a}=\bar{a}_k^+\rest w^\ast_0$
and $\bar{a}'_k=\bar{a}_k^+\rest w_\ast^1$.
Hence, $\langle \bar{a}\,\hat{}\,\bar{a}'_k:\,k\in {\Bbf Z}\rangle$ is an
indiscernible sequence,
and $\Rang(\bar{a})\cap\Rang(\bar{a}'_k)=\emptyset$ for all $k$. In addition,
and for $k_1\neq k_2$ we have $\Rang(\bar{a}'_{k_1})\cap
\Rang(\bar{a}'_{k_2})=\emptyset$ and $\Rang(\bar{a}\,\hat{}\,\bar{a}''_k)=X_k$.

Now we define a model $N$. Its universe is $\cup_{0\le l<3}
\{\cl_M(\bar{a}\,\hat{}\,
\bar{a}'_l\,\hat{}\,\bar{a}'_{l+1})\}$, and $Q_i^N=Q_i^M\cap N$, $F_j^N=
\cup\{F_{j,l}:\,l<3\}$, where
$F_{j,l}=F_j^M\rest \cl_M(\bar{a}\,\hat{}\,
\bar{a}'_l\,\hat{}\,\bar{a}'_{l+1})$, or $F_{j,l}=F_j^M\rest
(\cl_M(\bar{a}\,\hat{}\,
\bar{a}'_l\,\hat{}\,\bar{a}'_{l+1}))^2$, as appropriate.
Note that $N$ is well defined, and that it is a model of $T_0$.
$N$ is not necessarily a model of $T_0^+$, as the function $F_1$ may be only
partial.
Notice that $X_l\subseteq N$ for $l\in [0, 3]$. We wish to define $N'$ like
$N$, but identifying $\bar{a}_0^+$ and $\bar{a}^+_3$ coordinatwise.
We shall now check that this will give a well defined model of $T_0$. Note
that by the proof of Observation \ref{insert} we have
\begin{equation*}
\begin{split}
N'=& \bigcup_{0\le l<3} X_l\cup\bigcup_{0\le l<3}\{F_1^N(c,d):\,c,d\in X_l\cup
X_{l+1}\\
&\&\,\,\{c,d\}\nsubseteq X_l\,\,\&\,\,
\{c,d\}\nsubseteq X_{l+1}\,\,\&\,\,F_1^N(c,d)\notin X_l\cup X_{l+1}\}.\\
\end{split}
\end{equation*}

The possible problem is that $F_i^{N'}$ might not be well defined, i.e.
there could perhaps be a case defined in two distinct ways. We verify that
this does not happen, by discussing various possibilities.

\underline{Case 1}. For some $b\in \Rang(\bar{a}^+_0)$, say $b=\bar{a}_0^+(t)$,
$b'=\bar{a}^+_3(t)$ and $j\in \{0,2,3\}$ we have $F_j(b)\neq F_j(b')$ after the
identification of  $\bar{a}^+_0$
with $\bar{a}^+_3$. As
$\bar{a}^+_k$'s are closed, we have $F_j(b)=\bar{a}^+_0(s)$
and $F_j(b')=\bar{a}^+_3(s')$ for some $s,s'$. By indiscernibility,
we have $s=s'$, hence the identification will make $F_j(b)=F_j(b')$.

\underline{Case 2}. For some $s,t$ we have that $F_1(\bar{a}^+_0(s), \bar{a}^+_0(t))$
and $F_1(a^+_3(s), a^+_3(t))$ are well defined, but not the same
after the identification of $\bar{a}^+_0$ and $\bar{a}^+_3$. This case cannot
happen, as can be seen similarly as in the Case 1.

\underline{Case 3}. For some $\tau(x,y)\in \{F_1(x,y), F_1(y,x)\}$ and
$d_1=\bar{a}_0^+(s), d_2=\bar{a}_3^+(s)$ and some $e\in N$ we have that
$\tau^N(d_1,e),\tau^N(d_2,e)$ are well defined but do not get identified when
$N'$ is defined.

By Case 2, we have that $e\notin \bar{a}$ and $s\notin w_0^\ast$. 
As $\tau(e,d_1)$ is well defined and $d_1\in X_0\setminus\bar{a}$,
necessarily $e\in \cl_M(X_0\cup X_1)$. Similarly, as $\tau(e,d_2)$ is well
defined and $d_2\in X_3\setminus\bar{a}$, we have 
$e\in \cl_M(X_2\cup X_3)$. But, as $F_1(e,d_l)$ is well defined, we have $e\in
Q_2\cup Q_0$. Hence $e\in \cl_M(X_0\cup X_1)\setminus Q_1\subseteq X_0\cup
X_1$ and similarly $e\in X_2\cup X_3$. But this implies $e\in\bar{a}$,
a contradiction.

As $M$ is a
model of $T_0$, $F_0^M$ is onto (Claim \ref{amalg}(1)).
Suppose $y\in Q_0^N$, then for some $l\in [0,3)$ we have that
$y\in \cl_M(X_l\cup X_{l+1})$, so by Observation \ref{insert}, we have
$y\in X_l\cup X_{l+1}$. As each $X_l$ is closed in $M$, by Claim \ref{amalg}(6)
each $X_l$ is a model of
$T_0^+$, so $y\in \Rang(F^M_0)$, hence $y\in \Rang(F_0^N)$ and
$y\in\Rang(F_0^{N'})$. We can similarly prove that $F_3^{N'}$ is onto, and
as each $X_l$ is a model of $T^+_0$ we have by Claim \ref{amalg}(1) that
$Q_0^{N'}$, $Q_1^{N'}$ and $Q_2^{N'}$ are all non-empty. By Claim
\ref{amalg}(2), $N'$ can be extended to a model of $T^+_0$. 

By the choice of $\varphi$ and the fact that $T^\ast$ is complete we have that
\[
T^\ast\models(\forall \bar{x}_0,\bar{x}_1,\bar{x}_2)\neg [\varphi(\bar{x}_0,
\bar{x}_1)\wedge \varphi(\bar{x}_1,
\bar{x}_2) \wedge \varphi(\bar{x}_2,
\bar{x}_0)].
\]
As $T^\ast$ is the model completion of $T^+_0$, in particular $T^\ast$ and
$T^+_0$ are cotheories, so we have that 
\[
T^\ast\models(\forall \bar{x}_0,\bar{x}_1,\bar{x}_2)\neg [\varphi(\bar{x}_0,
\bar{x}_1)\wedge \varphi(\bar{x}_1,
\bar{x}_2) \wedge \varphi(\bar{x}_2,
\bar{x}_0)],
\]
yet in $N'$ we have
\[
N'\models\varphi(\bar{a}_0,
\bar{a}_1)\wedge \varphi(\bar{a}_1,
\bar{a}_2) \wedge \varphi(\bar{a}_2,
\bar{a}_0),
\]
by the identification of $\bar{a}_0$ and $\bar{a}_3$. This is a contradiction.
$\eop_{\ref{nsop}}$
\end{Proof of the Claim}

\begin{Definition}\label{Prs}
\begin{description}
\item{(1)} A theory $T$ is said to {\em satisfy}
the {\em oak property as exhibited by
a formula} $\varphi(\bar{z},\bar{y}, \bar{x})$ iff
for any $\lambda,\kappa$ there are $\bar{b}_\eta
(\eta\in{}^{\kappa>}\lambda)$ and $\bar{c}_\nu (\nu \in
{}^{\kappa}\lambda)$ and $\bar{a}_i (i<\kappa)$
such that
\begin{description}
\item{(a)} $[\eta\initial\nu\,\,\&\,\,\nu \in {}^\kappa\lambda]
\implies\varphi[\bar{c}_\nu,\bar{b}_\eta, \bar{a}_{\llg(\eta)}]$,
\item{(b)} If 
$\eta\in {}^{\kappa>}\lambda$ and $\eta\,\hat{}\,\langle\alpha\rangle\initial
\nu_1\in {}^\kappa\lambda$ and $\eta\,\hat{}\,\langle\beta\rangle\initial
\nu_2\in {}^\kappa\lambda$, while $\alpha\neq\beta$ and $i>\llg(\eta)$,
\underline{then} $\neg\exists \bar{y}\,
[\varphi(\bar{c}_{\nu_1}, \bar{y}, \bar{a}_i)
\wedge \varphi(\bar{c}_{\nu_2}, \bar{y}, \bar{a}_i)]$,
\end{description}
and in addition $\varphi$ satisfies
\begin{description}
\item{(c) $\varphi(\bar{z},\bar{y}_1, \bar{x})\wedge
\varphi(\bar{z},\bar{y}_2, \bar{x})\implies \bar{y}_1=\bar{y}_2$.
}
\end{description}
We allow for the replacement of ${\mathfrak C}_T$ by ${\mathfrak C}^{\rm eq}_T$
(i.e. allow $\bar{y}$ to be a definable equivalence class).

\item{(2)} We say that oak holds for $T$
if this is true for
some $\varphi$.
\end{description}
\end{Definition}

\begin{Observation} If some $\lambda,\kappa$ exemplify that
oak$(\varphi)$ holds, then so do all $\lambda,\kappa$.
(This holds by the compactness theorem).
\end{Observation}

\begin{Remark}
We shall not need to use this, but let us remark that
witnesses $\bar{a},\bar{b},\bar{c}$ to oak$(\varphi)$ can
be chosen to be indiscernible along an appropriate index set (a tree).
This can be proved using the technique as in \cite{Sh -c},
Chapter VII, which uses the compactness argument and an appropriate
partition theorem.
\end{Remark}

\begin{Claim}\label{nasa} $T^\ast$ has oak. 
\end{Claim}

\begin{Proof of the Claim} Let
\[
\varphi(z,y,x)\deq
Q_2(z)\wedge Q_1(y)\wedge Q_0(x) \wedge F_0(y)=x
\wedge F_1(z,x)=y.
\]
Clearly, (c) of Definition \ref{Prs} (1) is satisfied.
Given $\lambda,\kappa $,
we shall define a model $N=N_{\lambda,\kappa}$ of $T_0^+$.
This will be a submodel of $\frak C={\mathfrak C}_{T^\ast}$ such that
its universe
consists of
$Q_0^N\deq\{a_i:\,i<\kappa\}$ with no repetitions,
$Q_1^N\deq\{b_\eta:\,\eta\in{}^{\kappa>}\lambda\}$ with no repetitions and
$Q_2^N\deq\{c_\nu:\,\nu\in{}^{\kappa}\lambda\}$ with no repetitions,
while $Q_0,Q_1, Q_2$ are pairwise disjoint.
We also require that the following
are satisfied in ${\frak C}={\frak C}_{T^\ast}$:

\[
F_0(b_\eta)=a_{\llg (\eta)},
F_1(c_{\nu}, a_i)=
b_{\nu\rest i}
\]
and that $N$ is closed under $F_2$ and $F_3$. That such a choice is possible
can be seen by writing the corresponding type and using the saturativity
of $\frak C$.

We can check that $N\models T_0^+$, and that $N$
is a submodel of ${\frak C}$ when understood as a model of $T^+_0$.
Clearly, (a) from Definition \ref{Prs}(1) is
satisfied
for $\varphi$ and $a_i, b_\eta, c_\nu$ in place of
$\bar{a}_i, \bar{b}_\eta, \bar{c}_\nu$ respectively. To see
(b), suppose
that $\eta,\alpha,\beta,\nu_1,\nu_2$ and $i$ are as there, but $d$ is such that
$\varphi(c_{\nu_1},d, a_i)\wedge \varphi(c_{\nu_2},d, a_i)$.
Hence $F_1(c_{\nu_1}, a_i)=F_1(c_{\nu_2}, a_i)$, so $\nu_1\rest i=
\nu_2\rest i$, a contradiction. This shows that $\varphi$ is a witness for
$T^\ast$ having oak.
$\eop_{\ref{nasa}}$
\end{Proof of the Claim}

Finally, a remark showing why this research continues \cite{Sh 457}.
The readers unfamiliar with $T_{\rm feq}^\ast$ can skip to the
next section without loss of generality. We use the notation for
$T^\ast_{\rm feq}$ which was used in \cite{DjSh 692}, while the
fact that this is equivalent to the notation in \cite{Sh 457} was 
explained in \cite{DjSh 692}.

\begin{Remark} After renaming, ${\frak C}_{T^\ast_{\rm feq}}^{\rm eq}$
and ${\frak C}_{T^\ast}^{\rm eq}$ are isomorphically embeddable
into each other. To see this suppose that
$M$ is a model of $T^\ast_{\rm feq}$. Let $A=\{x_\alpha:\,\alpha<\alpha^\ast\}$
be a set of representatives of $E^M$-equivalence classes. By the construction
of $T^\ast_{\rm feq}$, for every finite $F\subseteq \alpha^\ast$, there is $z$
such that $\wedge_{\alpha\in F} F(x_\alpha,z)=x_\alpha$. By the saturativity
of ${\frak C}_{T^\ast_{\rm feq}}$, there is $z_A\in {\frak C}_{T^\ast_{\rm feq}}$
such that $\wedge_{\alpha <\alpha^\ast} F(x_\alpha,z_A)=x_\alpha$. By the
axioms of $T_{\rm feq}$ it follows that $A\neq A'\implies z_A\neq z_{A'}$.

Now
we define a model
$N=N_0[M]$ of $T^\ast$. Its universe is
\begin{equation*}
\begin{split}
\card{M}\cup (P^M/E^M)\cup\{z_A:\,&A \mbox{ is a set of
representatives}\\ 
&\mbox{of $E^M$-equivalence classes}\}.\\
\end{split}
\end{equation*}
We let $Q_0=P^M/E^M$, $Q_1^N=P^M$ and 
\[
Q_2^N=Q^M\cup \{z_A:\,A \mbox{ is a set of
representatives of $E^M$-equivalence classes}\}.
\]
The functions of $N$ are
defined as follows. We firstly let
$F_0(x)\deq x/E^M$ and $F_1(z,x/E)=F(x,z)$. Notice that $F_1$ is well defined,
because if $z\in M$ then certainly $F(x,z)\in M$ for all $x\in Q_1^N$
and if $z=z_A$ for some $A$ then the definition of $A$ guarantees that 
$F(x,z_A)\in M$ for every $x\in M$. Also, we have 
\[
F_0(F_1(z,x/E^M))=F_0(F(x,z))=x/E^M.
\]
It remains to define $F_2$ and $F_3$. Let us first see that $\card{Q_2}\ge
\card{Q_0}$. By the definition of $T^\ast_{\rm feq}$, each equivalence 
class of $M$ is infinite. Hence, the number of distinct
sets of representatives of the $E^M$-equivalence classes is at least
$\card{Q_0}^{\aleph_0}\ge \card{Q_0}$, and by the definition of $Q_2$ we have
$\card{Q_2}\ge\card{Q_0}$. We can choose $F_3$ as any onto function from
$Q_2$ to $Q_0$, and apply Claim \ref{amalg}(2). Hence $N$ is a model of
$T^+_0$ and can be seen as a submodel of ${\frak C}^{\rm eq}_{T^\ast}$.

Conversely, given $M$ a model of $T^\ast$, we define $N=N_1[M]$ by letting
its universe be $Q_1^M\,\bigcup Q^M_2$ and $P^N=Q_1^M$, while $Q^N=Q_2^M$.
We let 
\[
y\,E z\mbox{ iff } F_0^M(y)=F_0^M(z)\mbox{ and }F^N(x,z)=F_1(z, F_0(x)).
\]
We also let $x\,R\,z\iff F^N(x,z)=x$. It is easily seen that $N\models T_{{\rm
feq}}$.

Using this equivalence and the fact that oak and $NSOP_3$ are preserved
up to isomorphism of ${\frak C}^{\rm eq}$, we obtain:

\begin{Corollary}
\label{feq} (1) $T_{\rm feq}^\ast$ has oak.

{\noindent (2)} $T_{\rm feq}^\ast$ has $NSOP_3$.
\end{Corollary}

Part (2) of Corollary \ref{feq} was stated without proof in \cite{Sh 500}.

\end{Remark}

\section{The theorems} In this section we present two general
theorems showing that under certain cardinal arithmetic assumptions
oak theories do not admit universal models.

\begin{Theorem}\label{main} Assume that
\begin{description}
\item{(1)} $\cf(\kappa)=\kappa<\mu<\mu^+<\lambda=\cf(\lambda)$,
\item{(2)} $\lambda<\mu^\kappa$,
\item{(3)} $\kappa\le\sigma\le\lambda$,
\item{(4)} There are families ${\cal P}_1\subseteq [\lambda]^\kappa$
and ${\cal P}_2\subseteq [\sigma]^\kappa$ such that
\begin{description}
\item{(i)} for every $g:\,\sigma\into\lambda$ there is $X\in \PP_2$ with
$\{g(i):\,i\in X\}\in \PP_1$,
\item{(ii)} $\card{\PP_1}<\mu^\kappa, \card{\PP_2}\le\lambda$,
\end{description}
\item{(5)} $T$ is a theory of size $<\lambda$ which satisfies
oak$(\varphi(\bar{z},\bar{y},\bar{x}))$.
\end{description}
\underline{Then}
\[
\mbox{univ}(T,\lambda)\ge \mu^\kappa.
\]
\end{Theorem}

\begin{Definition} For cardinals $\kappa,\mu$ we define
\[
{\cal U}_{J^{\rm bd}_\kappa}(\mu)\deq\min\{\card{{\cal P}}:
\PP\subseteq [\mu]^\kappa\,\,\&\,\,(\forall b\in [\mu]^\kappa)(\exists
a\in \PP)(\card{a\cap b}=\kappa\}.
\]
\end{Definition}
\begin{Theorem}\label{mainprime} Assume that
\begin{description}
\item{(1)} $\cf(\kappa)=\kappa<\mu<\mu^+<\lambda=\cf(\lambda)$,
\item{(2)} $\lambda<{\cal U}_{J^{\rm bd}_\kappa}(\mu)$,
\item{(3)} $\kappa\le\sigma\le\lambda$,
\item{(4)} There are families ${\cal P}_1\subseteq [\lambda]^\kappa$
and ${\cal P}_2\subseteq [\sigma]^\kappa$ such that
\begin{description}
\item{(i)} for every $g:\,\sigma\into\lambda$ there is $X\in
{\cal P}_2$ such that $\card{
\{g(i):\,i\in X\}\cap Y}=\kappa$ for some
$Y\in {\cal P}_1$,
\item{(ii)} $\card{\PP_1}<{\cal U}_{J^{\rm bd}_\kappa}(\mu)$,
$\card{\PP_2}\le\lambda$, \end{description}
\item{(5)} $T$ has oak$(\varphi(\bar{z},\bar{y},\bar{x}))$.
\end{description}
\underline{Then}
\[
\mbox{univ}(T,\lambda)\ge {\cal U}_{J^{\rm bd}_\kappa}(\mu).
\]
\end{Theorem}

\begin{Proof} We shall use the same proof for both Theorem \ref{main}
and Theorem \ref{mainprime}. The two main Lemmas are the same for
both theorems, and we shall indicate the differences which occur
toward the end of the proof. Let $a_i\,(i<\kappa)$, $b_\eta\,(\eta\in
{}^{\kappa>}\lambda)$ and $c_\nu(\nu\in {}^\kappa\lambda)$ exemplify the
oak property of $\varphi(z,y,x)$ for $\lambda$ and $\kappa$.

For notational
simplicity, let us assume that $\llg
(\bar{x})=\llg(\bar{y})= \llg(\bar{z})=1$. Let
$\bar{C}=\langle C_\delta:\,\delta\in S\rangle$ for some $S\subseteq
S^\lambda_\kappa$
with $\otp{C_\delta}=\mu$
and $C_\delta$ a closed subset of
$\delta$,
be a club guessing sequence,
(i.e. for every $E$ a club of $\lambda$, there is
$\delta\in S$ with $C_\delta\subseteq E$) such that 
\[
\alpha<\lambda\implies\card{\{C_\delta\cap\alpha:\,\delta\in
S,\,\&\,\,\alpha\in\nacc(C_\delta)\}}<\lambda.
\]
Such a sequence exists by S. Shelah's \cite{Sh 420} (section \S1). For each $\delta$, let
$\langle \alpha_{\delta,\zeta}:\,\zeta<\mu\rangle$ be the increasing
enumeration of $C_\delta$. Let ${\frak C}^+$ be a
(saturated enough) expansion of ${\frak C}_T$ by
Skolem functions for ${\frak C}_T$.

\begin{Definition} (1) For $\bar{N}=\langle N_\gamma:\,\gamma\le\lambda
\rangle$ an $\elementary$-increasing continuous sequence of models of $T$
of size $\le\lambda$, and for $c,a\in \bigcup_{\gamma<\lambda}N_\gamma$,
and $\delta\in S$,
we let
\begin{equation*}
\begin{split}
{\rm inv}_{\bar{N}}(c, C_\delta, a)\deq
\{\zeta<\mu:\,&[(\exists b\in N_{\alpha_{\delta,\zeta+1}})(
N_\lambda\models\varphi[c,b,a])\,\,\&\,\,\\
&\neg
(\exists b\in
 N_{\alpha_{\delta,\zeta}})(
N_\lambda\models\varphi[c,b,a])\}.\\
\end{split}
\end{equation*}

 {\noindent (2)} For a set $A$ and $\delta$, $\bar{N}$ as above, let
\[
{\rm inv}^A_{\bar{N}}(c, C_\delta)\deq\bigcup\{ {\rm inv}_{\bar{N}}(c,
C_\delta, a):\,a\in A\}.
\]
\end{Definition}

\begin{Note} Notice that
${\rm inv}_{\bar{N}}(c, C_\delta, a)$ is always a singleton or empty
and that ${\rm inv}_{\bar N}^A(c,C_\delta)\in [\mu]^{\le\card{A}}$.
\end{Note}

\begin{Construction Lemma}\label{prva} For every
unbounded $A^\ast\in[\mu]^\kappa$ of order type $\kappa$,
there is an $\elementary$-increasing continuous sequence
$\bar{N}_{A^\ast}=\langle N^{A^\ast}_\gamma:\,\gamma<\lambda\rangle$
of models of $T$ of size $<\lambda$ and a set
$\{\hat{a}_i:\,i<\sigma\}$ of elements of $N^{A^\ast}_0$ such that for every
$X\in  \PP_2$, for every
$\delta\in S$ with $\min(C_\delta)$ large enough, 
there is $c\in N_{A^\ast}\deq\bigcup_{\gamma<
\lambda}N^{A^\ast}_\gamma$ such that
${\rm inv}{}^{\{\hat{a}_i:\,i\in X\}}_{\bar{N}_{A^\ast}}(c,C_\delta)=
A^\ast$. \end{Construction Lemma}

\begin{Proof of the Lemma} Let $\PP_2=\{X_\alpha:\,\alpha<\alpha^\ast\le
\lambda\}$.

Given $A^\ast$.
Let $f=f_{A^\ast}$ be the increasing
enumeration of $A^\ast$, so $f:\,\kappa\into\mu$. For $\delta\in S$ let
$\nu_\delta\deq\langle \alpha_{\delta,\zeta}:\,\zeta\in A^\ast\rangle$ be an
increasing enumeration,
hence $c_{\nu_\delta}$ is well defined, as is $b_\eta$ for $\eta\initial
\nu_\delta$. For $X\in \PP_2$, let $\rho_X$
be an increasing surjection from the successor ordinals
$<\kappa$ onto $X$. By a
compactness argument, we can see that there are
$\langle \hat{a}_i:\,i<\sigma\rangle$ and for $X\in \PP_2$, sequences
$\langle c^X_{\nu_\delta}:\,\delta\in S\rangle$,
$\langle b^X_\eta:\,\eta\initial\nu_\delta
\,\,\&\,\,\llg(\eta)\mbox{ a successor}\,\,\&\,\,\delta\in S\rangle$
such that
\[
\eta\initial \nu_\delta\implies \models\varphi[c^X_{\nu_\delta},
b^X_\eta,\hat{a}_{\rho_X(\llg(\eta))}]
\]
and the appropriate translation of (b) from Definition \ref{Prs} holds.
Let for $\gamma<\lambda$ the model $N_\gamma^{A^\ast}$ be the reduction to
${\cal L}(T)$ of the Skolem hull in ${\frak C}^+$ of
\[
\{\hat{a}_i:\,i\in \cup_{\alpha<\gamma}X_\alpha\}
\cup\{c^{X_\alpha}_{\nu_\delta}:\,\alpha<\gamma
\,\,\&\,\,\delta\in S\cap\gamma \,\,\&\,\,\sup(\Rang(
\nu_\delta))<\gamma\}\]
\[
\cup\{b_\eta^{X_\alpha}:\,\alpha<\gamma\,\,\&\,\,\eta\initial\nu_\delta
\mbox{ for some }\delta
\in S\cap\gamma\,\,\&\,\,\sup(\Rang(\eta))<\gamma
\,\,\&\,\,\llg(\eta)\mbox{ a successor}\}. \]
Hence $\langle \bar{N}^{A^\ast}=N^{A^\ast}_\gamma:\,\gamma<\lambda\rangle$ is
$\elementary$-increasing
continuous and for $\gamma<\lambda$ we have $\card{N^{A^\ast}_\gamma}<\lambda$.
The latter is true because in the last clause
\[
\card{\{C_\delta\cap\alpha:\,\delta\in
S,\alpha<\gamma\,\alpha\in\nacc(C_\delta)\}}<\lambda
\]
by the choice of
$\bar{C}$.
Given $\alpha<\alpha^\ast$, $X=X_\alpha$ and $\delta\in S$
with $\min(C_\delta)\ge\alpha+1$ we shall
show that with \[
I\deq{\rm inv}^{\{\hat{a}_i:\,i\in X\}}_{\bar{N}^{A^\ast}}(c^X_{\nu_\delta},
C_\delta)
\]
we have $I= A^\ast$.
Notice that
that $\varepsilon<\kappa\implies \alpha_{\delta, f(\varepsilon)}
>\alpha$.
Let $i\in X$
and let $\eta=\langle \alpha_{\delta, f(\varepsilon)}:\,\varepsilon
\le\beta\rangle$, where $\beta+1=\rho_X^{-1}(i)$. We have that 
$\eta\initial \nu_\delta$ and $i=\rho_X(\llg(\eta))$. Hence
$\varphi[c^X_{\nu_\delta}, b^X_\eta,  \hat{a}_i]$ holds. 
Let $\zeta=f(\beta)$. We then
have that
$b_\eta^X\in N^{A^\ast}_{\alpha_{\delta,\zeta}+1}\subseteq
N^{A^\ast}_{\alpha_{\delta,\zeta+1}}$
(as $\alpha_{\delta,\zeta+1}>\gamma$), but $b^X_\eta\notin
N^{A^\ast}_{\alpha_{\delta,\zeta}}$.
It follows from the property (c) of
Definition \ref{Prs} that $b^X_\eta$ is the only $b$ for which 
$\models\varphi[c^X_{\nu_\delta}, b, \hat{a}_i]$. Hence $\zeta
=f(\beta)\in I$. So $A^\ast\subseteq I$ because every element of $A^\ast$
is $f(\beta)$ for some $\beta$ as above.

In the other direction, suppose $\zeta\in I$ and let $i\in X$ be such that
$\zeta$ is in ${\rm inv}_{\bar{N}}(c^X_{\nu_\delta}, C_\delta, \hat{a}_i)$.
Hence for some $b\in N_{\alpha_\delta,\zeta+1}^{A^\ast}\setminus
N_{\alpha_\delta,\zeta}^{A^\ast}$ we have
$\models\varphi[c^X_{\nu_\delta}, b, \hat{a}_i]$. Constructing $\eta$ as in the
previous paragraph, we have
$\models\varphi[c^X_{\nu_\delta}, b^X_\eta, \hat{a}_i]$. In conclusion, using
property (c) of
Definition \ref{Prs} again, we see that
$b=B^x_\eta$ so $\zeta=f(\beta)$ for some $\beta$.
So $A^\ast=I$.

$\eop_{\ref{prva}}$
\end{Proof of the Lemma}

\begin{Preservation Lemma}\label{druga} Suppose that $N$ and $N^\ast$ are models of $T$ both
with universe $\lambda$, and $f:\,N\into N^\ast$ is an 
elementary embedding, while
$\langle N_\gamma:\,\gamma<\lambda\rangle$ and
$\langle N_\gamma^\ast:\,\gamma<\lambda\rangle$ are continuous increasing
sequences of models of $T$
of cardinality $<\lambda$ with $\bigcup_{\gamma<\lambda}N_\gamma=N$ and
$\bigcup_{\gamma<\lambda}N_\gamma^\ast=N^\ast$. Further suppose that $\{
\hat{a}_\alpha:\,\alpha<\kappa\}\subseteq N$ is given. Let
\[
E\deq\{\gamma:\,(N,N^\ast,f)\rest\gamma\elementary(N,N^\ast,f)\,\,\&\,\,
\sup(\{f(\hat{a}_\alpha):\,\alpha<\kappa\})<\gamma\},
\]
hence a club of $\lambda$.

\underline{Then} for every $c\in N$ and $\delta$ with $C_\delta
\subseteq E$ we have
\[
{\rm inv}^{\{\hat{a}_\alpha:\,\alpha<\kappa\}}_{\bar{N}}(c, C_\delta)=
{\rm inv}^{\{f(\hat{a}_\alpha):\,\alpha<\kappa\}}_{\bar{N}}(f(c), C_\delta).
\]
\end{Preservation Lemma}

\begin{Proof of the Lemma} Fix $c\in N$ and $\delta\in S$ as required, and let
$a=\hat{a}_\alpha$ for some $\alpha<\kappa$. We shall see that
${\rm inv}_{\bar{N}}(c, C_\delta,a)={\rm inv}_{\bar{N}^{\ast}}(f(c),
C_\delta, f(a))$.

Suppose $\zeta<\mu$ is an element of ${\rm inv}_{\bar{N}}(c, C_\delta,a)$, so
there is $b\in N_{\alpha_{\delta,\zeta+1}}$ with $N\models\varphi[c,b,a]$, while
there is no such $b\in N_{\alpha_{\delta,\zeta}}$. We have that
$N^\ast$ satisfies $\varphi[f(c),f(b), f(a)]$. As
$C_\delta\subseteq E$ we have that $\alpha_{\delta,\zeta+1}
\in E$, and as $b\in N_{\alpha_{\delta,\zeta+1}}$, clearly $f(b)\in
N^\ast_{\alpha_{\delta,\zeta+1}}$. Similarly, by the definition of $E$
again, we have $f(b)\notin N^\ast_{\alpha_{\delta,\zeta}}$.
By the assumptions on $\varphi$ we have
\[
N^\ast\models``(\forall y)[\varphi(f(c),y,f(a))\implies y=f(b)]",
\]
so
$\zeta\in
{\rm inv}_{\bar{N}^{\ast}}(f(c),
C_\delta, f(a))$.

In the other direction, suppose $\zeta<\mu$ is an element of
${\rm inv}_{\bar{N}^{\ast}}(f(c),
C_\delta, f(a))$, so there is $b^\ast\in N_{\alpha_\delta,\zeta+1}^\ast$
with $N^\ast\models\varphi[f(c),b^\ast,f(a)]$, while
there is no such $b^\ast\in N_{\alpha_\delta,\zeta}^\ast$. Hence
$N^\ast\models \exists y\,( \varphi[f(c),y,f(a)])$,
so $N\models\exists y\, (\varphi[c,y,a])$. Let $b\in
N$ be such that
$N\models \varphi[c,b,a]$. Hence $N^\ast\models \varphi[f(c), f(b), f(a)]$.
Again by (c) of
Definition \ref{Prs},
we have $f(b)=b^\ast$, so $b\in N_{\alpha_\delta,\zeta+1}\setminus
N_{\alpha_\delta,\zeta}$ by elementarity. As this $b$ is unique (by (c) of
Definition \ref{Prs}), we have
$\zeta\in {\rm inv}_{\bar{N}}(c, C_\delta,a)$.
$\eop_{\ref{druga}}$
\end{Proof of the Lemma}

{\em Proof of the Theorems continued.s} Theorem \ref{main}
[Theorem \ref{mainprime}]. To conclude the proof of the theorems,
given $\theta<\mu^\kappa$ [$\theta<{\cal U}_{J^{\rm bd}_\kappa}(\mu)$],
we shall see that
${\rm univ}(T,\lambda)>\theta$. Without loss of generality,
we can
assume that $\theta
\ge\lambda+\card{{\cal P}_1}$.
Given $\langle N^\ast_j:\,j<\theta\rangle$ a sequence of models of $T$ each
of size $\lambda$, we show that
these models are not jointly universal. So suppose they were.
Without loss of generality, the universe of each $N^\ast_j$ is $\lambda$. Let
$\bar{N}^\ast_j=\langle N_{\gamma,j}^\ast:\,\gamma<\lambda\rangle$ be an
increasing continuous sequence of models of $T$ of size $<\lambda$
such that $N^\ast_j
=\bigcup_{\gamma<\lambda} N^\ast_{\gamma,j}$,
for $j<\theta$. For each $A\in \PP_1$ (so $A\in [\lambda]^\kappa$),
$\delta\in S$, $j<\theta$ and $d\in N^\ast_j$, we compute ${\rm
inv}^{A}_{\bar{N}^\ast_j}(d,C_\delta)$, each time obtaining an element of
$[\mu]^{\le\kappa}$. The number of elements of $[\mu]^{\le\kappa}$ obtained in
this way is
\[
\le\card{{\cal P}_1}\cdot\card{S}\cdot\theta\cdot\lambda\le\theta.
\]
By the choice of $\theta$ [and the definition of ${\cal U}_{J^{\rm bd}_\kappa}
(\mu)$],
we can choose $A^\ast\in [\mu]^\kappa$ such that
$A^\ast$ is
not equal to any of these sets [is almost disjoint (i.e. has
intersection of size $<\kappa$) to all these
sets]. Let $N\deq N_{A^\ast}$ be as guaranteed to exist by the Construction
Lemma, and let
$\{\hat{a}_i:\,i<\sigma\}$ and
$\bar{N}_{A^\ast}\deq\langle N^{A^\ast}_\gamma:\,\gamma\le\lambda\rangle$ be as
in that Lemma.
Without loss of generality, by taking an isomorphic copy if necessary, the
universe of $N$ is $\lambda$.
Suppose that $j<\theta$ and $f:\,N\into N^\ast_j$ is an embedding,
and let \[
E\deq\{\delta<\lambda:\,(N,N^\ast_j,f)\rest\delta\elementary
(N,N^\ast_j,f)\}. \]
Let $g:\,\sigma\into\lambda$
be given by $g(i)=f(\hat{a}_i)$. Let $X=X_\alpha\in \PP_2$ be such that
$\{f(\hat{a}_i):\,i\in X\}\in \PP_1$,
[for some $Y\in {\cal P}_1$ we have $\card{\{f(\hat{a}_i):\,i\in X\}\cap
Y}=\kappa$], and let $c\in N$ be such that ${\rm inv}_{\bar{N}}^{\{
\hat{a}_i:\,i\in X\}}(c, C_\delta)=A^\ast$. By the Preservation Lemma, we have
${\rm inv}_{\bar{N}^\ast_j}^{\{f(\hat{a}_i):\,i\in X\}}(f(c),
C_\delta)= A^\ast$ [${\rm inv}_{\bar{N}^\ast_j}^{\{f(
\hat{a}_i):\,i\in X\}} (f(c), C_\delta)\cap A^\ast$ includes ${\rm
inv}_{\bar{N}^\ast_j}^{\{ f(\hat{a}_i):\,i\in X\}\cap Y}(f(c), C_\delta)\cap
a^\ast$, which has cardinality $\kappa$]. This is a contradiction with the
choice of $A^\ast$.

$\eop_{\ref{main}}$
\end{Proof}

\begin{Remark}\label{sigma} We comment on the assumptions used in Theorems
\ref{main} and \ref{mainprime}. Although the theorems do not use the
assumption $\cf(\mu)=\kappa$, this situation is the natural one for the
assumptions given. If $\cf(\mu)\le\kappa<\mu$ we have
${\rm pp}_{J^{\rm bd}_\kappa}(\mu)\le {\cal U}_{J^{\rm bd}_\kappa}(\mu)$. For
example, we have the following

\begin{Corollary}\label{posljedica} Let $T$ be a theory with the oak property.
Suppose that $\cf(\mu)=\kappa<\mu<\mu^+ <
\lambda=\cf(\lambda)$ and $\lambda<{\cal U}_{J^{\rm bd}_\kappa}(\mu)$
(e.g. pp${}_{J_\kappa^{\rm bd}}(\mu)>\lambda$) while $2^\kappa\le\lambda$, and
\begin{equation*}
\mbox{for some }n, {\rm cov}(\lambda,\kappa^{+n+1},\kappa^{+n+1},
\kappa^{+n})=\lambda\tag{$\ast_{\lambda,\kappa}$}\label{taggy}
\end{equation*}
\underline{then}
${\rm univ}(T,\lambda)\ge{\cal U}_{J^{\rm bd}_\kappa}(\mu)$.
\end{Corollary}

\begin{Proof} We use Theorem \ref{mainprime} with $\sigma=\kappa^{+n+1}$
for $n$ as in (\ref{taggy}). By the choice of $n$, there are ${\cal P}_1,
{\cal P}_2$ as
required and of cardinality $\lambda$.$\eop_{\ref{posljedica}}$
\end{Proof}

Note that the consistency of the failure of $(\ast_{\lambda,\kappa})$ for
any $\lambda\ge\kappa^{+\omega}, \kappa=\cf(\kappa)$ is not known, and that
for our purposes even weaker statements suffice. See \cite{Sh 460}.

If $\aleph_0<\kappa=\cf(\mu)$ and for all $\theta<\mu$
we have $\theta^\kappa<\mu$, then
\[
\mbox{pp}{}_{J^{\rm bd}_\kappa}(\mu)=
\mu^\kappa={\cal U}_{J^{\rm bd}_\kappa}(\mu)
\]
(by \cite{Sh -g}, Chapter VII, \S1).

If $\lambda>\kappa=\cf(\kappa)$ and $\sigma=\lambda$, if we cannot find 
${\cal P}_1$ and ${\cal P}_2$ as in Theorem \ref{mainprime}(i) with
$\card{{\cal P}_1}+\card{{\cal P}_2}\le\lambda$, \underline{then} for every
${\cal P}\subseteq [\lambda]^\kappa$ with $\card{\PP}\le\lambda$, we can find
$X\in [\lambda]^\lambda$ such that $(\forall a\in {\cal P})(\card{a\cap X}
<\kappa)$, which is a rather strong requirement.

Another comment is the
necessity of introducing
the
cardinal $\sigma$ at the outset of Theorem \ref{main}
and Theorem \ref{mainprime}.
In most instances of cardinal arithmetic, assuming that the
other requirements are satisfied, requirement (4) cannot be fulfilled
with $\kappa=\sigma$.
But if for example
$\lambda=\lambda^{[\sigma]}$ (for a definition see
\cite{Sh 460}; the equality holds e.g. if $\lambda<\aleph_{\sigma}$),
and $\kappa<\sigma$ is such that $\sigma^\kappa<\mu$, and for some
sequence $\langle\lambda_i:\,i<\kappa\rangle$ of regulars increasing to
$\mu$ the reduced product $(\Pi_{i<\kappa} \lambda_i/J^{\rm bd}_\kappa)$
is $\lambda^+$-directed, the  
assumptions of Theorem \ref{main} will hold. In fact, by
\cite{Sh 460} we have

\begin{Corollary} In Theorem \ref{main}, if (1)+(2) hold, while
$\kappa=\cf(\mu)<\beth_\omega\le\mu$, \underline{then} for every
large enough $\sigma\in (\kappa,\beth_\omega)$, parts (3) and (4)
of the assumptions of Theorem \ref{main} hold as well,
so univ$(T,\lambda)\ge \mu^\kappa$.
\end{Corollary}
\end{Remark}

\eject

\end{document}